\documentclass[11pt,reqno]{amsart}
\usepackage{amsmath,amssymb,amsthm,upref,graphicx,mathrsfs}

\usepackage{color,xcolor}
\usepackage[
  colorlinks=true,
  linkcolor=blue,
  citecolor=blue,
  urlcolor=blue]{hyperref}

\newtheorem{thm}{Theorem}[section]
\newtheorem{cor}[thm]{Corollary}
\newtheorem{lem}[thm]{Lemma}

\numberwithin{equation}{section}

\begin{document}

\leftline{ \scriptsize}

\vspace{1.3 cm}
\title
{determinants preserving maps on the spaces of symmetric matrices and skew-symmetric matrices}
\author{Ratsiri Sanguanwong $^{a}$ and Kijti Rodtes$^{a,\ast}$ }
\thanks{{\scriptsize
		$^{\ast}$ Corresponding author.\\
		E-mail addresses: ratsiris60@nu.ac.th (Ratsiri Sanguanwong), kijtir@nu.ac.th (Kijti Rodtes).\\
		$^{A}$ Department of Mathematics, Faculty of Science, Naresuan University, Phitsanulok 65000, Thailand.\\}}
\hskip -0.4 true cm

\maketitle


\begin{abstract}
Denote $\Sigma_n$ and $Q_n$ the set of all $n \times n$ symmetric and skew-symmetric matrices over a field $\mathbb{F}$, respectively, where $char(\mathbb{F})\neq 2$ and $\lvert \mathbb{F} \rvert \geq n^2+1$. A characterization of $\phi,\psi:\Sigma_n \rightarrow \Sigma_n$, for which at least one of them is surjective, satisfying
$$\det(\phi(x)+\psi(y))=\det(x+y)\qquad(x,y\in \Sigma_n)$$ is given. Furthermore, if $n$ is even and $\phi,\psi:Q_n \rightarrow Q_n$, for which $\psi$ is surjective and $\psi(0)=0$, satisfy
$$\det(\phi(x)+\psi(y))=\det(x+y)\qquad(x,y\in Q_n),$$
then $\phi=\psi$ and $\psi$ must be a bijective linear map preserving the determinant.
\end{abstract}

\vskip 0.2 true cm


\pagestyle{myheadings}
\markboth{\rightline {\scriptsize Ratsiri Sanguanwong and  Kijti Rodtes}}
{\leftline{\scriptsize }}
\bigskip
\bigskip


\vskip 0.4 true cm

\section{introduction}

The determinant preserving problem has been studying for a long time. The first characterization was investigated by Frobenius \cite{Frobenius1897} in 1897. The result stated that if $\phi: M_n(\mathbb{C})\rightarrow\mathbb{C}$ is a surjective linear transformation satisfying
\begin{equation}\label{assumption}
\det(a)=\det(\phi(a))\qquad (a\in M_n(\mathbb{C})),
\end{equation}
then it must be in the form
\begin{equation}\label{conclusion}
\phi(a)=uav \quad (a\in M_n(\mathbb{C})) \mbox{\ \ \ or\ \ \ }\phi(a)=ua^tv \quad (a\in M_n(\mathbb{C})),
\end{equation} 
where $u,v\in M_n(\mathbb{C})$ with $\det(uv)=1$.

In 2002, Dolinar and Šemrl \cite{DolinarSemrl2002} generalized Frobenius's result by dropping the linearity of the map and change (\ref{assumption}) into the form
\begin{equation}\label{assumption2}
\det(a+\alpha b)=\det(\phi(a)+\alpha\phi(b))\qquad(a,b\in M_n(\mathbb{C}),\alpha\in\mathbb{C}).
\end{equation}
However, the same conclusion as (\ref{conclusion}) is obtained.

This result was improved next year by Tan and Wang \cite{TanWang2003}. In that work, the result holds on arbitrary fields that satisfy some conditions of cardinality and characteristic. They proved that even if the surjective property is omitted, the assumption (\ref{assumption2}) also implies (\ref{conclusion}). In 2004, Cao and Tang \cite{CaoTang2004} published a similar work on the space of all symmetric matrices. The result on positive definite matrices can be found in \cite{Huangetal2016}.

In 2019, Costara accomplished another improvement. Assume that $\mathbb{F}$ is a field with $\lvert \mathbb{F}\rvert \geq n^2+1$ and $\phi,\psi:M_n(\mathbb{F})\rightarrow M_n(\mathbb{F})$ such that one of them is surjective and
\begin{equation}\label{assumption3}
\det(a+b)=\det(\phi(a)+\psi(b))\qquad (a,b \in M_n(\mathbb{F})).
\end{equation}
He proved that there must exists $x_0,u,v\in M_n(\mathbb{F})$ with $\det(uv)=1$ such that
$$\phi(x)=u(x+x_0)v\quad(x\in M_n(\mathbb{F})) \mbox{\ \ and \ \ }\psi(x)=u(x-x_0)v\quad(x\in M_n(\mathbb{F}))$$ or $$\phi(x)=u(x+x_0)^tv\quad(x\in M_n(\mathbb{F})) \mbox{\ \ and \ \ }\psi(x)=u(x-x_0)^tv\quad(x\in M_n(\mathbb{F})).$$

In this paper, denote $\Sigma_n = \Sigma_n(\mathbb{F})$ and $M_n = M_n(\mathbb{F})$. By using the idea of \cite{Costara2019}, the following theorem holds.

\begin{thm}\label{thm}
	Suppose that $\mathbb{F}$ is a field such that $char(\mathbb{F})\neq2$ and $\lvert\mathbb{F}\rvert \geq n^2+1$.
	Let $\phi,\psi:\Sigma_n \rightarrow \Sigma_n$, for which at least one of them is surjective, be satisfying
	\begin{equation}\label{assumption4}
	\det(\phi(x)+\psi(y))=\det(x+y)\qquad(x,y\in \Sigma_n).
	\end{equation}
	Then there exists $x_0\in \Sigma_n$, $u\in M_n$ and $\beta\in\mathbb{F}$ with $\det(\beta u^2)=1$ such that $$\phi(x)=\beta u(x+x_0)u^t\qquad(x\in \Sigma_n)$$
	and $$ \psi(x)=\beta u(x-x_0)u^t\qquad(x\in \Sigma_n).$$
\end{thm}

The proof of Theorem \ref{thm} is placed in the next section. Now, let $\alpha \in \mathbb{F}\setminus\{0\}$. Assume that $\phi,\gamma:\Sigma_n \rightarrow \Sigma_n$, for which at least one of them is surjective, satisfying
$$\det(\phi(x)+\alpha\gamma(y))=\det(x+\alpha y)\qquad(x,y\in \Sigma_n).$$ The equation (\ref{assumption4}) will be obtained if we define $\psi:\Sigma_n\rightarrow \Sigma_n$ by $\psi(x)=\alpha\gamma(\alpha^{-1}x)$. By Theorem \ref{thm}, it is not difficult to see that there exists $x_0\in \Sigma_n, u\in M_n$ and $\beta\in\mathbb{F}$ with $\det(\beta u^2)=1$ such that $$\phi(x)=\beta u(x+x_0)u^t\qquad(a\in \Sigma_n) \mbox{\ \ and \ \ } \gamma(x)=\beta u(x-\alpha^{-1}x_0)u^t\qquad(a\in \Sigma_n).$$
The following corollaries hold.
\begin{cor}
	Let $\phi,\psi:\Sigma_n \rightarrow \Sigma_n$, for which at least one of them is surjective, be satisfying
	\begin{equation*}
	\det(\phi(x)-\psi(y))=\det(x-y)\qquad(x,y\in \Sigma_n).
	\end{equation*}
	Then $\phi=\psi$ and there exists $x_0\in \Sigma_n,u\in M_n$ and $\beta\in\mathbb{F}$ with $\det(\beta u^2)=1$ with $\det(\beta u^2)=1$ such that $$\phi(x)=\beta u(x+x_0)u^t\qquad(x\in \Sigma_n).$$
\end{cor}
\begin{cor}
	Let $\alpha_1,\alpha_2$ be two distinct nonzero elements of $\mathbb{F}$ and $\phi,\psi:\Sigma_n \rightarrow \Sigma_n$, for which at least one of them is surjective, satisfying
	\begin{equation*}
	\det(\phi(x)+\alpha_i\psi(y))=\det(x+\alpha_i y)\qquad(x,y\in \Sigma_n, i=1,2).
	\end{equation*}
	Then $\phi=\psi$ and there exists $u\in M_n$ and $\beta\in\mathbb{F}$ with $\det(\beta u^2)=1$ such that $$\phi(x)=\beta uxu^t\qquad(x\in \Sigma_n).$$
\end{cor}

Let $Q_n$ denote the space of all skew-symmetric matrices over a field $\mathbb{F}$. For $n\in\mathbb{N}$, if $n$ is odd, all $n \times n$ skew-symmetric matrices are singular. This means that every mapping from $Q_n$ into itself is already preserving the determinant. Moreover, when $n$ is even, the following theorem is true.
\begin{thm}\label{thm2}
	Suppose that $\mathbb{F}$ is a field such that $char(\mathbb{F})\neq2$ and $\lvert\mathbb{F}\rvert \geq n^2+1$.
	Let $n$ be an even positive integer and $\phi,\psi:Q_n \rightarrow Q_n$, for which $\psi$ is surjective and $\psi(0)=0$, be satisfying
	\begin{equation*}\label{assumption5}
	\det(\phi(x)+\psi(y))=\det(x+y)\qquad(x,y\in Q_n).
	\end{equation*}
	Then $\phi=\psi$ and $\psi$ is a bijective linear map preserving the determinant.
\end{thm}

\section{Proof of main results}

Through out this section, the field $\mathbb{F}$ means any field with $char(\mathbb{F})\neq2$ and $\lvert\mathbb{F}\rvert \geq n^2+1$. According to the corollary of Theorem 6.35 \cite{FriedbergInselSpenceBook}, if $a$ is a symmetric matrix of rank $r$ over a field of characteristic not equal to 2, $a$ is congruent to a diagonal matrix having $r$ nonzero entries. This fact is used in the proof of the following lemma.
\begin{lem}\label{lem1}
	Let $a,b \in \Sigma_n$ be such that $\det(a+x)=\det(b+x)$ fot all $x \in \Sigma_n$. Then $a=b$. 
\end{lem}
\begin{proof}
	Denote $c=b-a$ and $y=a+x$. Then $\det(y)=\det(c+y)$ for all $y \in \Sigma_n$. Suppose that $\operatorname{rank}(c)=r$. There is an invertible matrix $p$ such that $$pcp^t=\operatorname{diag}(d_1,d_2,\dots,d_r,0,0\dots,0)=:d,$$ where $d_i\neq 0$ for all $i=1,\dots,r$. Let $w=\operatorname{diag}(0,0,\dots,0,w_1,\dots,w_{n-r})$ be a diagonal matrix having $n-r$ nonzero entries and denote $z=p^{-1}w(p^{-1})^t$. It is obvious that $d+w$ is invertible, and hence so is $c+z$. Because $\det(z)=\det(c+z)\neq 0$, we also have that $\det(w)\neq 0$. It follows that $r=0$. Thus $c=0$, which implies $a=b$. 
\end{proof}

In Lemma 2.1\cite{Costara2019}, every $n\times n$ matrix is represented by a row and a column vector of length $n^2$. Those vectors is used in order to find the trace of a product of any two matrices. Since the dimension of $M_n$ is $n^2$, an invertible $n^2 \times n^2$ matrix can be constructed by those vectors. We shall prove the following lemma by using that idea. Due to the dimension of $\Sigma_n$, the length of our vectors should be $\frac{n^2+n}{2}$. So, the constructions are a bit different from those in \cite{Costara2019}.
\begin{lem}\label{lem2}
	Let $\psi:\Sigma_n\rightarrow \Sigma_n$ and $\chi$ be a map defined on the set of all invertible symmetric matrices values into $\Sigma_n$. If $tr(xy)=tr(\chi(x)\psi(y))$ for all invertible symmetric matrices $x$ and all symmetric matrices $y$, then $\psi$ is linear.
\end{lem}
\begin{proof} For $a=(a_{ij})$ and $1\leq i,j \leq n$, denote
\begin{equation*}
c_{ij} = \left\{
\begin{array}{ll}
1, & \hbox{if $i=j$,} \\
\sqrt{2}, & \hbox{otherwise,}
\end{array}
\right.
\end{equation*}	
$$S_i(a)=(c_{ii}a_{ii}, c_{i(i+1)}a_{i(i+1)},\dots, c_{in}a_{in}),$$ and $$R_a=(S_1(a), S_2(a), \dots, S_n(a)).$$ We also denote $C_a = R_a^t$. Then, for any $a,b\in \Sigma_n$, 
\begin{eqnarray*}
R_aC_b 
	& = & \sum_{i=1}^n a_{ii}b_{ii}+2\sum_{i=1}^n\sum_{j=i+1}^n a_{ij}b_{ij}\\
	& = & \sum_{i=1}^n a_{ii}b_{ii}+\sum_{i=1}^n\sum_{j=i+1}^n a_{ij}b_{ji}+\sum_{i=1}^n\sum_{j=i+1}^n a_{ji}b_{ij}\\
	& = & \sum_{i=1}^n\sum_{j=1}^n a_{ij}b_{ji}\\
	& = & tr(ab).
\end{eqnarray*}
By the assumption, for each invertible symmetric matrix $a$ and each $b\in \Sigma_n$, we have that
\begin{equation}\label{eq1}
R_aC_b=R_{\chi(a)}C_{\psi(b)}.
\end{equation} 

Remark that the basis containing only invertible matrices for $\Sigma_n$ exists. To prove this fact, suppose that $B=\{a_1,\dots,a_k\}$ is a linearly independent subset of $\Sigma_n$ such that $a_1,\dots,a_k$ are invertible and every invertible symmetric matrix is contained in $\operatorname{Span}(B)$. Let $a$ be singular and $p(\alpha)=\det(a+\alpha e)$, where $e$ is the identity matrix. Because $p(\alpha)$ is a polynomial of degree $n$, it has at most $n$ distinct roots. Since $\lvert \mathbb{F} \rvert \geq n^2+1$, there is $\lambda \in \mathbb{F}$ such that $\det(a+\lambda e)=p(\lambda)\neq 0$. This implies that $a+\lambda e \in \operatorname{Span}(B)$. Note that $e \in \operatorname{Span}(B)$. This concludes that $a = (a+\lambda e) - \lambda e \in \operatorname{Span}(B)$. Thus, $\operatorname{Span}(B)=\Sigma_n$ which implies that $B$ is a basis.
 
Let $\{x_1,x_2,\dots,x_m\}$ be a basis for $\Sigma_n$ such that each $x_j$ is invertible. Define $$X=\begin{pmatrix}
R_{x_1}\\R_{x_2}\\\vdots\\R_{x_m}.
\end{pmatrix} \mbox{ and } Z=\begin{pmatrix}
R_{\chi(x_1)}\\R_{\chi(x_2)}\\\vdots\\R_{\chi(x_m)}
\end{pmatrix}.$$
Then $XC_y = ZC_{\psi(y)}$ for all $y \in \Sigma_n$. Now, let $\{y_1,y_2,\dots,y_m\}$ be a basis for $\Sigma_n$. 

Putting $$Y=(C_{y_1} \ C_{y_2} \ \dots \ C_{y_m}) \mbox{ and } W=(C_{\psi(y_1)} \ C_{\psi(y_2)} \ \dots \ C_{\psi(y_m)}).$$ By (\ref{eq1}), $ZW=XY$. Because $X$ and $Y$ are invertible, $Z$ is also invertible. Let $Q=Z^{-1}X$. We can conclude that $C_{\psi(y)}=QC_y$ for all $y \in \Sigma_n$. Therefore $\psi$ is linear.
\end{proof}

For every $x \in \Sigma_n$, the $(i,j)$-cofactor of $x$ is denoted by $A_{ij}(x)$. For each $1\leq i,j \leq n$, denote $e_{ij}$ the matrix which $(i,j)$-entry is $1$ and other entries are $0$. It is obvious that $A_{ii}(a+te_{ii}) = A_{ii}(a)$ for all $t \in \mathbb{F}$ and all $1 \leq i \leq n$. If $i \neq j$, $A_{ji}(a+e_{ij})$ can be expressed as the following lemma.

\begin{lem}\label{lem3} 
	Let $a\in M_n$ and $t\in\mathbb{F}\setminus\{0\}$. If $i\neq j$, then $A_{ji}(a+te_{ij})=A_{ji}(a)-t\det(a^{\{ij\}})$, where $a^{\{ij\}}$ denote the submatrix of $a$ obtained by deleting row $i,j$ and column $i,j$. 
\end{lem}
\begin{proof}
	For each $x\in M_n$, denote $x(i\vert j)$ a submatrix of $x$ obtained by deleting row $i$ and column $j$. Note that, for each $x\in M_n$ and $t\in\mathbb{F}$, $\det(x+te_{ij})=\det(x)+tA_{ij}(x)$. Then
	\begin{eqnarray*}
		A_{ji}(a+te_{ij})
		 & = & (-1)^{j+i}\det((a+te_{ij})(j\vert i))\\
		 & = & (-1)^{j+i}\det(a(j\vert i)+te_{ij}(j\vert i))
	\end{eqnarray*}
	Since $i \neq j$, we have that the nonzreo entry of the $(n-1) \times (n-1)$ matrix $te_{ij}(j\vert i)$ must be $(i,j-1)$-entry if $i<j$ and it is $(i-1,j)$-entry if $i>j$. This means that
	\begin{eqnarray*}
		A_{ji}(a+te_{ij})
		& = & (-1)^{j+i}\det(a(j\vert i)+te_{ij}(j\vert i))\\
		& = & (-1)^{j+i}(\det(a(j\vert i))+t(-1)^{i+j-1}\det(a^{\{ij\}}))\\
		& = & A_{ji}(a)-t\det(a^{\{ij\}}).
	\end{eqnarray*}
Thus, the proof is completed.
\end{proof}

By using these lemmas, the proof of Theorem \ref{thm} can be accomplished.
\\

\textbf{Proof of Theorem \ref{thm}}
	This theorem is trivial for $n=1$. Suppose that $n\geq 2$. Without loss of generality, we may suppose that $\psi$ is surjective. First, we will prove the case $\psi(0)=0$. By the assumption, taking $y=0$ yields $\det(\phi(x))=\det(x)$ for all $x\in \Sigma_n$. By the same argument as in the proof of Theorem 1.1 \cite{Costara2019}, we have $A_{ii}(\phi(\alpha a))$ is a polynomial of degree at most $n-1$, in the form
	\begin{equation*}\label{eq1101}
		A_{ii}(\phi(\alpha a)) = \det(\alpha a+ b_{ii}) - \alpha^n\det(a),
	\end{equation*}
	where $b_{ii}\in \Sigma_n$ such that $\psi(b_{ii})=e_{ii}$. For each $1\leq i <j \leq n$, let $b_{ij}\in\Sigma_n$ such that $\psi(b_{ij})=e_{ij}+e_{ji}$. Then
	\begin{eqnarray*}
	\det(\alpha a+ b_{ij}) & = & \det(\phi(\alpha a)+\psi(b_{ij}))\\
							& = & \det(\phi(\alpha a)+e_{ij}+e_{ji})\\
							& = & \det(\phi(\alpha a)+e_{ij}) +A_{ji}(\phi(\alpha a)+e_{ij})\\
							& = & \det(\phi(\alpha a))+A_{ij}(\phi(\alpha a))+A_{ji}(\phi(\alpha a)+e_{ij})\\
							& = & \det(\alpha a)+A_{ij}(\phi(\alpha a))+A_{ji}(\phi(\alpha a)+e_{ij})\\
							& = & \alpha^n\det(a)+A_{ij}(\phi(\alpha a))+A_{ji}(\phi(\alpha a)+e_{ij}),
	\end{eqnarray*}
	for all $a\in\Sigma_n$ and all $\alpha\in \mathbb{F}$.
	By Lemma \ref{lem3},
	$$\det(\alpha a+ b_{ij})=\alpha^n\det(a)+2A_{ij}(\phi(\alpha a))-\det((\phi(\alpha a))^{\{ij\}}).$$
	Let $d_{ij}\in \Sigma_n$ such that $\psi(d_{ij})=-(e_{ij}+e_{ji})$. So,
	\begin{eqnarray*}
	\det(\alpha a+ d_{ij})
	& = & \alpha^n\det(a)-A_{ij}(\phi(\alpha a))-A_{ji}(\phi(\alpha a)-e_{ij})\\
	& = & \alpha^n\det(a)-A_{ij}(\phi(\alpha a))-(A_{ij}\phi(\alpha a)+\det((\phi(\alpha a))^{\{ij\}})\\
	& = & \alpha^n\det(a)-2A_{ij}(\phi(\alpha a))-\det((\phi(\alpha a))^{\{ij\}}).
	\end{eqnarray*}
By combining these equations,
\begin{equation*}\label{eq1102}
A_{ij}(\phi(\alpha a)) = \frac{1}{4}(\det(\alpha a+ b_{ij}) - \det(\alpha a+d_{ij}))
\end{equation*}
is a polynomial with respect to $\alpha$ of degree at most $n-1$. Since our field can be finite and $char(\mathbb{F}) \neq 2$, the notation $\frac{1}{4}$ means $(1+1+1+1)^{-1}$.
For each $1\leq i,j,k \leq n$, let $\beta_{i,j,k-1}(a)$ be the coefficient of $\alpha^{k-1}$ of $A_{ij}(\phi(\alpha a))$, that is,
$$A_{ij}(\phi(\alpha a)) = \beta_{i,j,0}(a)+\alpha\beta_{i,j,1}(a)+\cdots+\alpha^{n-1}\beta_{i,j,n-1}(a).$$
Let $a$ be an invertible symmetric matrix. Then, for each nonzero $\alpha \in \mathbb{F}\setminus\{0\}$, we have
\begin{eqnarray*}
(\phi(\alpha^{-1} a))^{-1} & = & \frac{1}{\det(\phi(\alpha^{-1} a))}(A_{ji}(\phi(\alpha^{-1}a)))_{ij}\\
			& = & \frac{\alpha^n}{\det(a)}(\omega_0(a)+\alpha^{-1}\omega_1(a)+\cdots+ \alpha^{-(n-1)}\omega_{n-1}(a)),
\end{eqnarray*}
where $\omega_k(a) = (\beta_{j,i,k}(a))_{ij}$. For each $1\leq i \leq n$, denote $\chi_i=\frac{1}{\det(a)}\omega_{n-i}$. We obtain that
$$(\phi(\alpha^{-1}a))^{-1}=\alpha\chi_1(a)+\alpha^2\chi_2(a)+\cdots+\alpha^n\chi_n(a).$$
For any $y \in \Sigma_n$,
\begin{eqnarray*}
\det(\alpha^{-1}a)\det(e+\alpha a^{-1}y) 
	& = & \det(\alpha^{-1}a + y)\\
	& = & \det(\phi(\alpha^{-1}a)+\psi(y))\\
	& = & \det(\phi(\alpha^{-1}a))\det(e+(\phi(\alpha^{-1}a))^{-1}\psi(y)),
\end{eqnarray*}
where $e$ is the identity matrix. Note that $\det(\alpha^{-1}a)=\det(\phi(\alpha^{-1}a))$. Then
\begin{eqnarray*}
\det(e+\alpha a^{-1}y)
	& = & \det(e+(\phi(\alpha^{-1}a))^{-1}\psi(y))\\
	& = & \det(e+(\alpha\chi_1(a)+\alpha^2\chi_2(a)+\cdots+\alpha^n\chi_n(a))\psi(y)).
\end{eqnarray*}
Consider the polynomial $$p(\alpha)=\det(e+(\alpha\chi_1(a)+\alpha^2\chi_2(a)+\cdots+\alpha^n\chi_n(a))\psi(y))-\det(e+\alpha a^{-1}y).$$
We can see that $p(\alpha)$ is a polynomial of degree at most $n^2$. However, each element of $\mathbb{F}$ is a root of $p(\alpha)$. Thus, it is a zero polynomial. This implies that $$\det(e+(\alpha\chi_1(a)+\alpha^2\chi_2(a)+\cdots+\alpha^n\chi_n(a))\psi(y))=\det(e+\alpha a^{-1}y).$$
Recall that for any matrix $b$, $tr(b)$ is equal to the coefficient of $\alpha$ of the polynomial $\det(e+\alpha b)$. Thus, $tr(a^{-1}y)=tr(\chi_1(a)\psi(y))$. Define a function $\chi$ on the set of all invertible symmetric matrices into $\Sigma_n$ by $\chi(a)=\chi_1(a^{-1})$. For each invertible symmetric matrix $a$ and each $y\in \Sigma_n$, we have  $$tr(ay)=tr(\chi(a)\psi(y)).$$
By Lemma \ref{lem2}, $\psi$ is linear.

Claim that $\phi(0)=0$. By the assumption, $\det(\phi(0)+\psi(y))=\det(y)$ for all $y \in \Sigma_n$. Let $x=\alpha^{-1}y$. We derive that
\begin{eqnarray*}
\alpha^n\det(x) 
	& = & \det(\alpha x)\\
	& = & \det(\phi(0)+\psi(\alpha x))\\
	& = & \det(\phi(0)+\alpha\psi(x))\\
	& = & \alpha^n\det(\alpha^{-1}\phi(0)+\psi(x)).
\end{eqnarray*}
This implies that, for each nonzero $\mu\in \mathbb{F}$, $\det(\mu\phi(0)+\psi(x))=\det(x)$. Thus the polynomial $\det(\alpha\phi(0)+\psi(x))-\det(x)$, which is of degree at most $n$, must be a zero polynomial. Consequently, $\det(\psi(x))=\det(x)$. For each $y \in \Sigma_n$, there exists $x\in \Sigma_n$ such that $\psi(x)=y$, and hence, $$\det(\phi(0)+y)=\det(x)=\det(\psi(x))=\det(y).$$
By Lemma \ref{lem1}, $\phi(0)=0$ as the claim. Now, for each $x,y \in \Sigma_n$ and $\alpha \in \mathbb{F}$, we have
$$\det(x+\alpha y)=\det(\psi(x+\alpha y))=\det(\psi(x)+\alpha\psi(y)).$$
Due to the result of Cao and Tang \cite{CaoTang2004}, there exists an invertible matrix $u$ and $\beta\in\mathbb{F}$ with $\det(\beta u^2)=1$ such that $\psi(a)=\beta uau^t$ for all $a \in \Sigma_n$.

To evaluate $\phi$, we first show that $\psi$ is a bijection. Let $a,b \in \Sigma_n$ such that $\psi(a)=\psi(b)$. For each $c \in \Sigma_n$,
$$\det(c+a)=\det(\phi(c)+\psi(a))=\det(\phi(c)+\psi(b))=\det(c+b).$$ By Lemma \ref{lem1}, $a=b$. Thus, $\psi$ is a bijective, which implies that $\psi^{-1}$ exists. Moreover, $\psi^{-1}$ is linear and preserving the determinant. So, for each $x,y \in \Sigma_n$,
\begin{eqnarray*}
	\det(x+y)
	& = & \det(\phi(x)+\psi(y))\\
	& = & \det(\psi^{-1}(\phi(x)+\psi(y)))\\
	& = & \det(\psi^{-1}\phi(x)+y).
\end{eqnarray*}
By using Lemma \ref{lem1} again, $\phi(x)=\psi(x)=\beta uxu^t$. The case $\psi(0)=0$ is accomplished. 

Suppose that $\psi(0)=a\neq 0$. Let $\psi_a = \psi-a$ and $\phi_a=\phi+a$. This yields that $\phi_a(0)=0$ and $$\det(\phi_a(x)+\psi_a(y))=\det(\phi(x)+\psi(y))=\det(x+y)\qquad(x,y\in \Sigma_n).$$ Consequently, there exists $u\in M_n$ and $\beta\in\mathbb{F}$ with $\det(\beta u^2)=1$ such that $$\phi(x)=\beta uxu^t-a \qquad(x\in \Sigma_n) \mbox{\ \ and \ \ } \psi(x)=\beta uxu^t+a\qquad(x\in \Sigma_n).$$
Putting $x_0=-\beta^{-1}u^{-1}a(u^t)^{-1}$. The proof is completed.

In order to prove Theorem \ref{thm2}, the skew-symmetric version of Lemma \ref{lem1} and Lemma \ref{lem2} are needed.

\begin{lem}\label{lem4}
	Let $n$ be an even positive integer and $a,b \in Q_n$ be such that $\det(a+x)=\det(b+x)$ fot all $x \in Q_n$. Then $a=b$. 
\end{lem}
\begin{proof}
	Denote $c=b-a$ and $y=a+x$. Then $\det(y)=\det(c+y)$ for all $y \in Q_n$. Note that the rank of skew-symmetric must be even. Suppose that $\operatorname{rank}(c)=2r$. It is known that two skew-symmetric are congruent if and only if they have the same rank (cf. \cite{Albert1938}). There must exist an invertible matrix $p$ such that $$pcp^t=\operatorname{diag}(\underbrace{x,x,\dots,x}_{r \mbox{\ times}},0,0\dots,0),$$ where $x = \begin{pmatrix}
	0 & 1 \\ -1 & 0
	\end{pmatrix}$. Let $$w=\operatorname{diag}(0,0,\dots,0,\underbrace{x,x,\dots,x}_{\frac{n-2r}{2} \mbox{\ times}})$$ and denote $z=p^{-1}w(p^{-1})^t$. By the same argument as the proof of Lemma \ref{lem1} $a=b$. 
\end{proof}

\begin{lem}\label{lem5}
	Let $n$ be even and $\psi:Q_n\rightarrow Q_n$ and $\chi$ be a map defined on the set of all invertible skew-symmetric matrices values into $Q_n$. If $tr(xy)=tr(\chi(x)\psi(y))$ for all invertible skew-symmetric matrices $x$ and all skew-symmetric matrices $y$, then $\psi$ is linear.
\end{lem}
\begin{proof} For $a=(a_{ij})$ and $1\leq i,j \leq n$, denote
	$$S_i(a)=(c_{i(i+1)}a_{i(i+1)},\dots, c_{in}a_{in}),$$ and $$R_a=(S_1(a), S_2(a), \dots, S_{n-1}(a)).$$ We also denote $C_a = -R_a^t$. Then, for any $a,b\in Q_n$, 
	\begin{eqnarray*}
		R_aC_b 
		& = & -\sum_{i=1}^n \sum_{j=i+1}^n a_{ij}b_{ij}\\
		& = & -\sum_{i=1}^n \sum_{j=i}^n a_{ij}b_{ij}\\
		& = & \frac{1}{2}\sum_{i=1}^n \sum_{j=i}^n a_{ij}(-b_{ij})+\frac{1}{2}\sum_{i=1}^n \sum_{j=i}^n (-a_{ij})b_{ij}\\
		& = & \frac{1}{2}\sum_{i=1}^n \sum_{j=i}^n a_{ij}b_{ji}+\frac{1}{2}\sum_{i=1}^n \sum_{j=i}^n a_{ji}b_{ij}\\
		& = & \frac{1}{2}\sum_{i=1}^n \sum_{j=1}^n a_{ij}b_{ji}\\
		& = & \frac{1}{2}tr(ab).
	\end{eqnarray*}
	By the assumption, for each invertible skew-symmetric matrix $a$ and each $b\in Q_n$, we have that
	\begin{equation}\label{eq2}
	R_aC_b=R_{\chi(a)}C_{\psi(b)}.
	\end{equation} 
	
	By the same argument as the proof of Lemma \ref{lem2}, but replacing $e$ by $\operatorname{diag}(x,x,\dots,x)$, where $x = \begin{pmatrix}
	0 & 1 \\ -1 & 0
	\end{pmatrix}$, the basis containing only invertible skew-symmetric matrices exists. Let $\{x_1,x_2,\dots,x_m\}$ be a basis for $Q_n$ such that each $x_j$ is invertible, where $m = dim(Q_n) = \frac{n^2-n}{2}$. Define $$X=\begin{pmatrix}
	R_{x_1}\\R_{x_2}\\\vdots\\R_{x_m}.
	\end{pmatrix} \mbox{ and } Z=\begin{pmatrix}
	R_{\chi(x_1)}\\R_{\chi(x_2)}\\\vdots\\R_{\chi(x_m)}
	\end{pmatrix}.$$
	Then $XC_y = ZC_{\psi(y)}$ for all $y \in Q_n$. Now, let $\{y_1,y_2,\dots,y_m\}$ be a basis for $Q_n$. 

	Putting $$Y=(C_{y_1} \ C_{y_2} \ \dots \ C_{y_m}) \mbox{ and } W=(C_{\psi(y_1)} \ C_{\psi(y_2)} \ \dots \ C_{\psi(y_m)}).$$ By (\ref{eq1}), $ZW=XY$. Because $X$ and $Y$ are invertible, $Z$ is also invertible. Let $Q=Z^{-1}X$. We can conclude that $C_{\psi(y)}=QC_y$ for all $y \in Q_n$. Therefore $\psi$ is linear.
\end{proof}

So far,  Theorem \ref{thm2} can be proven.

\textbf{Proof of Theorem \ref{thm2}}
	By the assumption, taking $y=0$ yields $\det(\phi(x))=\det(x)$ for all $x\in Q_n$.
For each $1\leq i<j \leq n$, denote $e_{ij}$ the matrix which $(i,j)$-entry is $1$ and other entries are $0$. Then $e_{ij}-e_{ji}\in Q_n$. Because $\psi$ is surjective, there exists $b_{ij},d_{ij}\in Q_n$ such that $\psi(b_{ij})=e_{ij}-e_{ji}$ and $\psi(d_{ij})=e_{ji}-e_{ij}$. 

By using the idea of the proof of Theorem \ref{thm} and the fact that $A_{ij}(x)=-A_{ji}(x)$ for all $x \in Q_n$, for each $\alpha\in \mathbb{F}$,
$$\det(\alpha a+ b_{ij})=\alpha^n\det(a)+2A_{ij}(\phi(\alpha a))+\det((\phi(\alpha a))^{\{ij\}})$$
and
$$\det(\alpha a+ d_{ij})=\alpha^n\det(a)-2A_{ij}(\phi(\alpha a))+\det((\phi(\alpha a))^{\{ij\}}),$$
%
%
%
and hence, we can write the polynomial $A_{ij}(\phi(\alpha a))$ in the form
\begin{eqnarray*}
	A_{ij}(\phi(\alpha a)) & = & \frac{1}{4}(\det(\alpha a+ b_{ij}) - \det(\alpha a+d_{ij}))\\
	& = & \beta_{i,j,0}(a)+\alpha\beta_{i,j,1}(a)+\cdots+\alpha^{n-1}\beta_{i,j,n-1}(a).
\end{eqnarray*}

Now, the remaining of the proof is exactly the same arguments as in the proof of Theorem \ref{thm}.

\section*{Acknowledgements}
The authors would like to thank anonymous referee(s) for reviewing this manuscript.


\begin{thebibliography}{20}
	
	\bibitem{Albert1938}
	A.A. Albert. \textit{Symmetric and Alternate Matrices in An Arbitrary Field, I}. Transactions of the American Mathematical Society. 43(3) (1938), 386-436.
	
	\bibitem{CaoTang2004}
	C. Cao and X. Tang. \textit{Determinant preserving transformations on symmetric matrix spaces}. Electronic Journal of Linear Algebra. 11 (2004), 205-211.
	
	\bibitem{Costara2019}
	C. Costara. \textit{Nonlinear determinant preserving maps on matrix algebras}. Linear Algebra and its Applications. 583 (2019), 165-170.
	
	\bibitem{DolinarSemrl2002}
	G. Dolinar and P. Šemrl. \textit{Determinant preserving maps on matrix algebras}. Linear Algebra and its Applications. 348 (2002), 189-192.
	
	\bibitem{FriedbergInselSpenceBook}
	S.H. Friedberg, A.J. Insel and L.E. Spence Johnson. \textit{Linear algebra}. (4th ed.). Pearson Education, Inc., 2003.
	
	\bibitem{Frobenius1897}
	G. Frobenius. \textit{Über die Darstellung der endlichen Gruppen durch lineare Substitutionen}. Sitzungsber.
	Deutsch. Akad. Wiss. (1897), 994-1015.
	
	\bibitem{HornJohnsonBook}
	R.A. Horn and C.R. Johnson. \textit{Matrix analysis}. (2nd ed.). Cambridge university press, 2013.
	
	\bibitem{Huangetal2016}
	H. Huang, C.N. Liu, P. Szokol, M.C. Tsai, and J. Zhang. \textit{Trace and determinant preserving maps of matrices}. Linear Algebra and its Applications. 507 (2016), 373-388.
	
	\bibitem{StollBook}
	R. Stoll. \textit{Linear algebra and matrix theory}. (Dover edition). Dover publication, Inc., New York, 1952.
	
	\bibitem{TanWang2003}
	V. Tan and F. Wang. \textit{Trace and determinant preserving maps of matrices}. Linear Algebra and its Applications. 369 (2003), 311-317.
\end{thebibliography}
\end{document}